\documentclass{article}
\begin{document}

\title{ Hilbert modular forms and the Ramanujan conjecture}

\author{ Don Blasius}

\maketitle

\newtheorem{intro}{Theorem}

\noindent Let F be a totally real field. In this paper we study
the Ramanujan Conjecture for Hilbert modular forms and  the
Weight-Monodromy Conjecture for the Shimura varieties attached to
quaternion algebras over $F$. As a consequence, we deduce, at all finite
places of the field of definition, the full automorphic
description conjectured by Langlands of the zeta functions of
these varieties.  Concerning the first problem, our main result is
the following:

\begin{intro} The Ramanujan conjecture holds at all finite places for any cuspidal 
holomorphic automorphic representation $\pi$ of $GL(2,\mathbf{A}_F)$ having weights 
all congruent modulo 2 and at least 2 at each infinite place of $F$.

\end{intro}

\noindent  See below (2.2) for a more precise statement. For background, we note that the 
above result has been known for any such $\pi$ at all but finitely many places, and without 
the congruence restriction, since 1984 (\cite{brylab}), as a consequence of the direct local 
computation of the trace of Frobenius on the intersection cohomology of a Hilbert modular 
variety. Additionally, the local method of \cite{hmfzeta} is easily seen to yield the  result at 
{\emph{all}} finite places, for the forms $\pi$ which satisfy the  restrictive hypothesis that 
either $[F:\mathbf{Q}]$ is odd or  the local component $\pi_v$ is discrete series at some 
finite place $v$. Hence, the novel cases in  Theorem 1  are essentially those of the forms 
$\pi$ attached to $F$ of even degree, and which belong to the principal series at all finite 
$v$. \\

\noindent To prove Theorem 1,  we here proceed globally, using the fact 
(\cite{hmfzeta}, \cite{ohta}, \cite{rthmf}, \cite{ordinary}) that there exist
two dimensional irreducible (\cite{brhmf}, \cite{rthmf2}) $l$-adic representations 
$\rho^T_l(\pi)$ of the Galois group of $\overline{F}$ over $F$ attached
to such forms $\pi$. Crucial to us is the fact that these representations satisfy the Global 
Langlands Correspondence, i.e. that  at every finite place $v$ whose residue characteristic is different from $l$, the representations of the 
Weil-Deligne group defined by $\pi_v$ and $\rho^T_l(\pi)$ (\cite{hmfzeta},\cite{rthmf}, 
\cite{brhmf}, \cite{rthmf2}, \cite{ordinary}) are isomorphic. Thus we get  information about $\pi_v$ from 
that about the local Galois representation
$\rho^T_l(\pi)|D_v$  whenever we realize $\rho^T_l(\pi)$, or a closely related
representation $\rho^{\prime}_l(\pi)$, in some $l$-adic cohomology. Many such realizations 
are provided by the Shimura varieties attached to inner forms of $GL(2)/F$, and to the 
unitary groups $GU(2)/K$  and $GU(3)/K$ where $K$ is a totally real solvable extension of 
$F$. Actually,   to go beyond the case of lowest discrete series at $\infty$,  in order to 
obtain cohomological realizations of these Galois representations $\rho^T_l(\pi)$ it is 
necessary   to consider fiber systems of abelian varieties over these  unitary Shimura 
varieties.
However, we need no explicit treatment of them here since the result is contained in 
\cite{brhmf}. The fact that these are realizations  of $\rho^T_l$ follows from suitable local  
Hasse-Weil zeta function computations at all but finitely many {\em good places}; it is important 
to note that in this paper no new such computations at  bad places are done. 

\noindent To actually get the results, there are  several overlapping methods:\\

\noindent A.  If one of the weights is greater than 2, or if either (a) $[F:\mathbf{Q}]$ is 
odd  or (b) there is a finite place at which $\pi$ is discrete series, the result follows easily 
from a basic theorem of De Jong (\cite{alter}), the Local Langlands Correspondence, and 
the classification of unitary representations of $GL(2)$ over a local field. In all these cases 
there is a direct realization of $\rho^T_l(\pi)$ as a subquotient of an $l$-adic cohomology 
group of a variety. \\

\noindent B. If all the weights are 2,  we proceed, using a known case of Langlands 
functoriality,  by finding a geometric realization of a Galois representation 
$\rho^{\prime}_l(\pi)$, made using $\rho^T_l(\pi)$, and from which we can deduce crucial 
constraints on the Frobenius eigenvalues of $\rho^T_l(\pi)$ at an unramified place under 
study. While several approaches are possible, we here use one for which the L-function of 
$\rho^{\prime}_l(\pi)$ is, after a formal base change to a field $L$, a Rankin product L-
function defined by $\pi|_L$ and a Galois twist $^{\tau}\pi|_L $.  Unlike case (A) above, 
to conclude Ramanujan by an extension of that method we use
a stronger, global  Ramanujan estimate (\cite{sharam}) for $GL(2)$ which the local 
analytic theory cannot provide. Although several alternative constructions of $\rho^{\prime}_l(\pi)$ are possible,  the present method has the merit that, further developed,  it enables progress on the $p$-adic 
analogue of the Langlands correspondence for these forms. Nevertheless, in order not to 
obscure the simple formal structure of the paper, we defer $p$-adic questions to a sequel. \\ 

\noindent C. If all the weights are 2, we can give (See 4.2) prove Theorem 1 by a geometric
 argument (found after that of B.) using the fact that the Weight-Monodromy Conjecture is a theorem for surfaces. We  give both arguments since,  the method of B., although a little longer,  has a chance to be applicable to other cases, such as regular algbraic forms on $GL(N)$ where $N>2$. \\

\noindent In this paper, we have restricted our study to the case of forms having weights all 
congruent modulo 2. However, the method may extend to all holomorphic forms whose 
weights are all at least 2 at the infinite places. A key fact, already present in \cite{brhmf}, is 
that a suitable twist $\pi^{\prime}=\pi_K\otimes \chi$ of a CM quadratic base change  
$\pi_K$ of $\pi$ defines motivic forms on  appropriate unitary groups $GU(2)$ and 
$GU(3)$.  Once the Global Langlands Correspondence (See below, Section 2.3)is known 
for these forms, the Ramanujan Conjecture will follow by the methods of this paper. One 
natural approach is to generalize, in the setting of those $GU(2)$ which define curves, the 
results of Carayol (\cite{hmfzeta}), and then to extend by congruences (\cite{rthmf}), to the 
general case.\\

\noindent  The second main goal of the paper is to  provide new examples, of arbitrary 
dimension, and with $N$ (See the text for definitions) of many different, often highly 
decomposable, types, of the Weight-Monodromy Conjecture (WMC) (\cite{dnice}).
\begin{intro}

Let $Sh_{B}$ be  the Shimura variety attached to a quaternion algebra $B$ over a totally 
real field $F$. Then WMC holds for the $l$-adic cohomology of $Sh_{B}$ at all finite places $v$ whose residue characteristic is different from $l$.

\end{intro}

\noindent {\bf {Remarks.}}\\

\noindent 1. $Sh_{B}$ is a projective limit of varieties $Sh_{B,W}$ , where $W$ is an open 
compact subgroup of the finite adele group of the reductive $\mathbf{Q}$-group 
$G=G_B=Res_{F/\mathbf{Q}}(B*)$ associated to the multiplicative group of $B$. Each $Sh_{B,W}$ is defined over the canonically  defined 
number field $F^{\prime}$,  named by Shimura the {\em{reflex field}}; the definition is recalled below in 
Section 3. We say that WMC holds for $Sh_{B}$ if it holds for each smooth variety 
$Sh_{B,W}$.  \\

\noindent 2. The Shimura variety is not proper exactly when  $B=M_2(F)$, in which case the 
connected components of the $Sh_{B, W}$ are the classical Hilbert modular varieties. In 
this case, the theorem is understood to refer to the $l$-adic intersection cohomologies of the 
Baily-Borel compactification of $Sh_{B, W}$.\\

\noindent 3. Several authors have recently made significant progress  on  cases of WMC involving Shimura varieties.  In \cite{itopadic},  instances of WMC are shown for certain Shimura varieties $Sh$ associated to unitary groups.  In fact,   WMC is shown at places $v$ at which  $Sh$ admit $p$-adic uniformization.   In \cite{deshalit},  the $p$-adic extension of WMC is shown for a similar class of varieties: here $v$ divides  $p$. As already noted, this is a case not treated at all in this paper.  Finally, in \cite {tywmc} Taylor  and Yoshida  establish WMC, by careful study of the Rapoport-Zink spectral sequence, for all Shimura varieties  associated to the unitary groups defined by division algebras over a CM field which are definite at all but one infinite place. This is the key class studied in \cite{ht}, and is a vast generalization  of the Shimura curves studied in \cite{hmfzeta}.  As a consequence,  WMC is true for the $l$-adic representations attached to the class of essentially self-dual regular automorphic cusp forms on $GL(N, \mathbf{A}_F)$.  This result implies Theorem 1 for $\pi$ which are discrete series at some finite place,  in which case the result is due to Carayol (\cite {hmfzeta}).\\

\noindent As a corollary of the above result, we achieve easily the third main goal of the 
paper: the proof of Langlands' conjecture (\cite{zssv}) which describes, in automorphic terms,  the Frobenius semisimplification of the action of a decomposition group $D_v$ for $v$ on the $l$-adic Galois  cohomology of the quaternionic Shimura varieties. Here $v$ is any finite place of  the reflex field, and $l$ is a prime different from the residue characteristic of $v$. . This result  completes the zeta function computations of Langlands (\cite{zssv}), Brylinksi-Labesse  (\cite{brylab}) and Reimann  (\cite{quat})). 

\begin{intro} Let $B$ be a quaternion algebra over a totally real field $F$ having $B_v \cong M_2(\mathbf{R})$ for $r>0$ infinite places $v$ of $F$. Let $F^{\prime}$ be the 
canonical field of definition of the r-fold $Sh_B$ attached to $B$. Let $\pi^{\prime}$ be a 
cuspidal holomorphic representation of $G)=(B\otimes \mathbf{A}_F)^{*}$ such that 

\begin{enumerate}

\item $\pi^{\prime}_v$ has weight 2 at each split infinite place,

\item $\pi^{\prime}_v$ is one-dimensional at each ramified infinite place,

\item the central character $\omega$ of $\pi^{\prime}$ has the form $\omega=|\cdot|^{-
1}\Psi$, with a character of finite order $\Psi$.

\end{enumerate}

\noindent Let $l$ be a rational prime. Then for each finite place $v$ of $F^{\prime}$ whose residue characteristic is different from $l$, the isomorphism class of the 
Frobenius semisimple parameter $(\rho^*_{W,v}, N_{W,v})$ of the Weil-Deligne  group 
$WD_{v}$ of $F^{\prime}_{v}$ defined by the restriction to a decomposition group for 
$v$ of the action of $Gal(\overline{\mathbf{Q}}/F^{\prime})$ on

$$H^{r}(Sh_{B, W}, \overline{\mathbf{Q}_l})(\pi^{\prime}_{f,W})$$

\noindent coincides with the class of

$$ m(\pi^{\prime}_f, 
W)r_B(\sigma(JL(\pi^{\prime})_{p})|_{WD_{F^{\prime}_{v}}}),$$

\noindent where $m(\pi^{\prime}_f, W)$ is defined by $$ dim(H^{r}(Sh_{B, W}, 
\overline{\mathbf{Q}_l})(\pi^{\prime}_{f,W}))= 2^{r} m(\pi^{\prime}_f, W).$$

\end{intro}

\noindent Here, for $p$ the place of $\mathbf{Q}$ lying under $v$, 

\begin{enumerate}

\item

$JL(\pi^{\prime})_{p}$ is the $p$-component of the cuspidal representation of $GL(2, 
\mathbf{A}_F)$, obtained from $\pi^{\prime}$ via the Jacquet-Langlands correspondence 
JL. 

\item

 $\sigma(JL(\pi^{\prime})_{p})$ is the homomorphism of  $WD_p$ into the L-group 
$^{L}G$ which is, as usual, identified with the L-group of the $\mathbf{Q}$-group $R_{F/\mathbf{Q}}(GL(2))$.

\item $r_B$ is the complex representation of dimension $2^r$ defined by Langlands.

\item Let $\cal{H}_{W}$ be the level $W$ Hecke algebra of $G$ which consists of the convolution algebra of left and right $W$ invariant compactly supported functions on $G(\mathbf{A}_f)$. Then $\pi^{\prime}_{f, W}$ is the representation of $\cal{H}_{W}$ on the subspace $(\pi^{\prime}_{f})^{W}$ of $\pi^{\prime}_f $ consisting of the vectors fixed by all of $W$.

\end{enumerate}
\noindent For an exposition of (2), see \cite{zeta}, 3.5, and \cite{kull}. For an exposition of 
(3), defined by Langlands (\cite{zssv}), see \cite{zeta} esp. 5.1, 7.2. Definitions are briefly 
recalled as needed in the paper. Note that we are computing the L-functions as Euler 
products over the primes of $F^{\prime}$, not as Euler products over primes of 
$\mathbf{Q}$.  \\

\noindent This result may  be the first verification, for some  Shimura varieties of dimension 
greater than one,  at all places and levels,  of Langlands' general conjecture. Nevertheless, for 
the last two theorems, the proofs are rather formal and do not involve new direct local 
verifications of difficult facts. On the contrary, one key principle is that the 
semisimplification of the {\em{global}}, i.e. $Gal(\overline{\mathbf{Q}}/F^{\prime})$,  
Galois action on the $l$-adic cohomology of any variety in Theorem 2 is computable in simple 
ways from globally  {\em{irreducible}} $l$-adic representations which satisfy WMC and the 
Global Langlands correspondence at each place. Of course this type of fact does not hold 
locally: the WMC concerns, for each place $v$ of $F^{\prime}$, the nature of the 
associated indecomposable, and in general non-irreducible, Frobenius semisimple 
representations of the Weil-Deligne group.\\

\noindent {\bf{Acknowledgements}} . I thank the SFB 478 at
M\"{u}nster and the Universit\'{e} Henri Pasteur at Strasbourg for
their generous hospitality. I thank Ron Livn\'{e} for asking me
several years ago about the RC at {\emph{all}} places,  and Takeshi
Saito for a stimulating conversation which led me to highlight WMC in the 
presentation of the results. Finally, I thank the referee for conscientious reading of the manuscript and the reporting of many misprints.

\section{Background}

\subsection{Weil Numbers.}

\noindent Let $q$ be power of a rational prime $p$. An  {\it {integral $q$-Weil number of 
weight j}} $ \in {\mathbf Z}$ is an algebraic integer $\alpha$ having the property that, for  
each automorphism  $\sigma$ of  $\overline{\mathbf{Q}}$, we have 
$$|\sigma(\alpha)|=q^{j/2},$$
\noindent with a fixed $j$ independent of $\sigma$. We omit reference to $q$ or the 
weight $j$ when convenient. An algebraic number of the form $\beta=\alpha q^n$, for 
some $n \in  \mathbf{Z}$  and an integral $q$-Weil Number  $\alpha$ is called  a {\it{$q$-Weil number}}, or simply a Weil number, if the q is clear from context. Obvious facts 
about Weil numbers include:  (i) the $q$-Weil numbers form a group under 
multiplication;(ii) if $q=q_0^f$, then $\alpha $ is a $q$-Weil number of weight $j$ if and 
only if  $\sqrt[f]{\alpha}$ is a  $q_0$-Weil number of weight $j$; (iii) all roots of unity  
are $q$-Weil numbers of weight $0$ for all $q$.

\subsection {$l$-adic Representations.}

\noindent  Let $ K$ be a field and let $\Gamma_K =Gal(\overline{K}/K) $ be the group of 
$K$-linear automorphisms of its algebraic closure $\overline{K}$, endowed with the usual 
topology. For a prime $l$, let  $V$ be a finite dimensional vector space over $\overline{\mathbf{Q}_l}$, and let $\rho :\Gamma_K \rightarrow GL(V) $  be a homomorphism. We 
say that $\rho$ is an {\it{$l$-adic representation}} if there exists a finite extension $T$ of  
$\mathbf{Q}_l$, a $T$ vector space $V_0$, and a continuous homomorphism $\rho_0 :\Gamma_K 
\rightarrow GL(V_0) $ which becomes isomorphic to $\rho$ after extension of scalars on 
$V_0$ from $T$ to  $\overline{\mathbf{Q}_l}$. We use the notation $\rho$,$V$, and 
$(V, \rho)$ at will to denote such a representation. An $l$-adic representation $(V, \rho)$ is 
called {\em{motivic}} if there is a smooth projective variety $X$ over $K$ such that $(V, 
\rho)$ is isomorphic to a Tate twist of a subquotient of the $\Gamma_K$-module 
$H^{*}(\overline{X}, \overline{\mathbf{Q}_l})$ where $\overline{X}$ is the scalar 
extension of $X$ to the algebraic closure $\overline{K}$ of $K$. Here, for a 
$\Gamma_K$-module $(V, \rho)$, and $m \in \mathbf{Z}$, the Tate twisted module is the 
pair $(V(m), \rho(m))$ where  $V(m)=V$,  $\rho(m)=\rho\otimes \chi_l^m$, and 
$\chi_l$  is the usual $l$-adic cyclotomic character.

\subsection{Local Weil group.}

\noindent For the rest of this paper,  $K$ will denote a local field of characteristic 0 and residue characteristic $p$. 
We denote by  $q$ the number of the residue field. Of course, $q$ is a power of $p$. We let $l$ be any rational prime different from $p$. We recall 
some basics about  the Weil group $W_K \subset \Gamma_K$ of $K$.   Let $I$ be the 
inertia subgroup of $W_K$. Then $W_K/I$ is  isomorphic to the subgroup 
$q^{\mathbf{Z}}$ of $\mathbf{Q}^{*}$; the isomorphism is that induced by the 
homomorphism that sends an element $w$ of $W_K$ to the power $|w|$ of $q$ to which 
it raises the prime-to-$p$ roots of unity in the maximal unramified extension of K. Any 
element $\Phi$ of $W_K$  for which $|\Phi|=q^{-1}$ is called a {\em{Frobenius}}. Let 
$I_w$ be the subgroup of wild inertia, i.e. the maximal pro-$p$ subgroup of $I$. Let $I_t$ 
denote the quotient $I/I_w $ and let $W_{K, t}$ denote $W_K/I_w$. We call these 
groups the {\em tame inertia} group and the {\em tame Weil} group,  respectively. Then 
$I_t$ is non-canonically isomorphic to the product

$$ \prod_{l\neq p}\mathbf{Z}_l, $$

\noindent and $W_{K,t}$ is isomorphic to the semidirect product of  $\mathbf{Z}$ and 
$I_t$ ; the action of $W_K$ on $I_t$  is given by $$wxw^{-1}=|w|x $$
\noindent  for all $x \in I_t$ and all $w\in W_K$.  Choose, once and for all,  an isomorphism $t=(t_l)_{l\neq p}$ of  $I_t$  with

$$ \prod_{l\neq p}\mathbf{Z}_l. $$

\noindent   Let $(V, \rho)$  be an $l$-adic  representation of $\Gamma_K$. We extend, 
replacing $\Gamma_K$ by $W_K$,  the definition of an $l$-adic representation to $W_K$, 
and  thus each $l$-adic representation $\rho$ of $\Gamma_K$ gives rise, by restriction, to an $l$-adic 
representation of $W_K$ which we also denote by $(V, \rho)$.

\subsection {Grothendieck's Theorem}

\noindent  According to a basic result of Grothendieck (\cite{serretate}, Appendix), there is 
a  subgroup  $J$ of finite index  in $I$ such that, for $\sigma \in J$,

$$\rho(\sigma)=exp(t_l(\sigma)N)$$

\noindent  where $N\in End(V)$ is a uniquely determined nilpotent endomorphism.\\

\noindent  If we can take $J=I$ in this theorem, $(V, \rho)$   is said to be 
{\em{semistable}}. It is well-known that there exists a finite extension $L$ of $K$ such 
that $(\rho|_L , V)$ is semistable.

\subsection{Weil-Deligne parametrization of $l$-adic representations}

\noindent  Fix a choice $\Phi$ of a Frobenius in $W_K$. Define, for this $\Phi$, and any 
$\sigma$ in $W_K$,  an automorphism $$\rho_{WD}(\sigma)= \rho(\sigma)exp(-t_l 
(\Phi^{-log_q(|\sigma|)}\sigma) N)$$ of $V$.
Then  $\sigma \rightarrow \rho_{WD}(\sigma)$ is a continuous representation of $W_K$ 
whose restriction to $I$ has finite image. The triple
$(V, \rho_{WD}, N)$  depends on the choice of  $t_l$ and $\Phi$. Such a triple 
$(V^{\prime}, \rho_{WD}^{\prime}, N^{\prime})$  arises from an $l$-adic representation 
on $V$ of $\Gamma_K$ if and only if the relation
$$ \rho_{WD}^{\prime}(\sigma)N^{\prime}\rho_{WD}^{\prime}(\sigma)^{-
1}=|\sigma|N^{\prime},$$
\noindent  holds for all $\sigma \in W_K$.
\noindent Note that  $(V, \rho)$ is semistable if and only if $ \rho_{WD}$ is unramified, 
i.e. trivial on $I$.

\subsection{Frobenius semisimplification.}

\noindent  Following Deligne ([D3], 8.5),  let 
$\rho_{WD}(\Phi)=\rho_{WD}(\Phi)^{ss}u$ be the Jordan decomposition of 
$\rho_{WD}(\Phi)$ as the product of a diagonalizable matrix   $\rho_{WD}(\Phi)^{ss}$ 
and a unipotent matrix $u$. Define, for $\sigma \in W_K$,  $$\rho_{WD}^{ss}(\sigma)= 
\rho_{WD}(\sigma)u^{log_q(|\sigma|)}.$$

\noindent  Then $\rho_{WD}^{ss}$ is a semisimple representation of $W_K$ and, for all 
$\sigma $, $\rho_{WD}^{ss}(\sigma)$ is semisimple. The representation $(V, 
\rho_{WD}^{ss})$  is called the $\Phi$-semisimplification of $(V, \rho_{WD})$ and the 
triple $(V, \rho_{WD}^{ss}, N)$ is called the $\Phi$-semisimplification of $(V, 
\rho_{WD}, N)$.\\

\noindent  Let now  $\iota_l$ be an isomorphism of  $\overline{\mathbf{Q}_l}$ with the 
complex numbers $\mathbf{C}$. This will be fixed in any discussion, and, to avoid a 
cumbersome notation, we will identify $\overline{\mathbf{Q}_l}$ with $\mathbf{C}$, 
suppressing explicit reference  to $\iota_l$. We will use $\iota_l$ to define complex 
representations of the Weil-Deligne group (c.f. \cite{dcef}, 8.3, \cite{ntb}, 4.1, or 
\cite{ecwd}), via the triples $(V, \rho_{WD}^{ss}, N)$.

\subsection{$WD_K$.} The {\em{Weil-Deligne}} group $WD_K$ of  $K$ is the 
semidirect product of $W_K$ with $\mathbf{C}$ defined by the relation

$$\sigma z {\sigma}^{-1}= |\sigma|z$$

\noindent  for all $\sigma \in W_K$ and $z \in \mathbf{C}$. Using $\iota_l$ we regard 
$V$ as a finite dimensional complex vector space (i.e. if  $z \in \mathbf{C}$ and $v \in 
V$, we put $zv=\iota_l^{-1}(z)v$). Then $\rho_{WD}^{ss}$ is a continuous 
representation of $W_K$ on $V$, and $N$ is a nilpotent endomorphism of V. The 
complex triple $(V, \rho_{WD}^{ss}, N)$ defines, in view of (1.5), a representation 
$\rho^{*}$ of $WD_K$ by the rule

$$  \rho^*((z, \sigma))= exp(z N)\rho_{WD}^{ss}(\sigma)$$

\noindent  for all $(z, \sigma) \in WD_K$. Then $\rho^{*}$ satisfies \\

\noindent (i) the restriction to $W_K$ is semisimple, and\\

\noindent (ii) the restriction to $\mathbf{C}= \mathbf{G}_a(\mathbf{C})$ is algebraic.\\

\noindent  We denote the family of all complex representations satisfying (i) and (ii) by 
$Rep^{s}(WD_K)$ and denote members by pairs $(V,\rho^{\prime}) $;  a triple giving 
rise to $(V, \rho^{\prime})$ by the construction above is given by

$$(V, \rho^{\prime}|_{W_K}, N_{\rho^{\prime}})$$

\noindent  where $log$ of a unipotent matrix $M$ is the standard polynomial in $M-1$ 
inverting exponentiation on nilpotents and

$$N_{\rho^{\prime}}=log(\rho^{\prime}((1, 1))).$$\\

\noindent  Henceforth an element of $Rep^s(WD_K)$ is identified with the triple it 
defines. Note that a member of  $Rep^{s}(WD_K)$  is actually a semisimple 
representation of $WD_K$ if and only if it  factors through the quotient $W_K$. A 
member of $Rep^{s}(WD_K)$  is called {\em{semistable}} if it is trivial on $I$. We 
denote by $Rep^{ss}(WD_K)$  the subfamily of $Rep^{s}(WD_K)$ consisting of 
semistable representations. Of course, if $(V, \rho)$ is a semistable $l$-adic representation if 
and only if the associated element $(V,\rho^{\prime})$ of $Rep^{s}(WD_K)$ 
belongs to $Rep^{ss}(WD_K)$. \\

\noindent  As shown in \cite{dcef}, the isomorphism class of the $(V,\rho^{\prime}) \in 
Rep^{s}(WD_K)$   gotten from an $l$-adic representation of $\Gamma_K$ is independent 
of the choices of $\Phi$ and $t_l$. The class of $(V,\rho^{\prime})$ does depend on the 
choice of $\iota_l$, but, since any $\iota_l^{\prime}$ has the form $\iota_l^{\prime}= \eta 
\iota_l$  for an automorphism $\eta$ of $\mathbf{C}$, we see that, after such a change, 
$\rho^{\prime} $  is just replaced by the conjugate $\eta \rho^{\prime}$.

\subsection{Structure of semistable modules.}

\noindent Recall that a $WD_K$-module is {\em{indecomposable}} if it cannot be written 
as the direct sum two proper submodules. We have the following basic structure results 
([Roh]) for the members of  $Rep^{ss}(WD_K)$:\\
 
\noindent (i) Any member of $Rep^{ss}(WD_K)$  is isomorphic to a direct sum of 
indecomposable modules, hence of $V_{\alpha, t}$'s. As such the decomposition is unique 
up to re-ordering the factors, and replacing factors by isomorphic factors.\\ 
\noindent (ii) Any indecomposable member of $Rep^{ss}(WD_K)$ is isomorphic to 
exactly one of the form $V_{\alpha, t}= (\mathbf{C}^{t+1}, \rho_{\alpha, t}, N_t) $, 
where $\alpha$ is a non-zero complex number,  $t$ is a non-negative integer, and 
$\rho_{\alpha, t}$ is the unramified representation of $W_K$ defined by the rule: 
$$\rho_{\alpha, t}(\Phi)=Diag(\alpha, q^{-1}\alpha, ..., q^{-t}\alpha), $$
where $Diag$ denote diagonal matrix, and $N=(n_{ij})$, where $n_{ij}=0$ unless 
$i=j+1$, in which case $n_{ij}=1$.

\subsection{ Structure of Frobenius semisimple modules.} We have:\\

\noindent (i) Any member of $Rep^{s}(WD_K)$  is a direct sum of indecomposable 
submodules. As such the decomposition is unique up to re-ordering the factors, and 
replacing factors by isomorphic factors.\\ 

\noindent(ii) Any indecomposable representation is isomorphic to one of the form 
$V_{\Lambda,t}\stackrel{def}{=}\Lambda\otimes V_{q^{t/2},t}$ where $\Lambda$ (and 
hence t) is a uniquely determined irreducible representation of $W_K$, and any such  
representation is indecomposable. Such a representation is irreducible iff $t=0$.\\ 

\noindent (iii)  if $\Lambda $ is an irreducible representation of $W_K$ and $\Phi$ is any 
Frobenius element in $W_K$, and $\alpha$ is an eigenvalue of $\Phi$ in $\Lambda$, then 
$|\alpha|$ is independent of $\alpha$. \\

\noindent To see the last claim, note that  we can find a Galois extension $L$ of $K$ such 
that the restriction to $WD_L \subseteq WD_K$ of $\Lambda$ is unramified, hence a 
direct sum of unramified characters $\chi_k$. Since $\Lambda$ is irreducible, the $\chi_k$  
are permuted transitively by the natural action of $\Gamma(L/K)$. Regarding them, via 
local class field theory,  as characters of $L^*$, and letting  $\tau$ be an element of 
$\Gamma(L/K)$, the action is just that sending $\chi_k$ to $\chi_k \circ \tau =\chi_k$. 
Hence all $\chi_k$ are the same character $\chi$. Now let $\chi_0$ be an unramified 
character of $W_K$ such that $\chi_0 \circ N_{L/K}=\chi$, and consider the irreducible 
representation $\Lambda_0=\Lambda \otimes \chi_0^{-1}$. Then the restriction to $L$ of 
$\Lambda_0$ is trivial, and hence $\Lambda_0$ has finite image. In particular, 
$\Lambda(\Phi)=\Lambda_0(\Phi)\chi_0(\Phi)$, and so each eigenvalue $\alpha$ of $\Phi$ 
in $\Lambda$ is of the form  $\alpha= \zeta \chi_0(\Phi)$ with a root of unity $\zeta$. This 
proves (iii).\\

\noindent Let $\Lambda$ be an irreducible representation of $WD_K$. We call the real number $ w(\Lambda)=2log_{q}(|\alpha|)$, where $\alpha$ is 
any eigenvalue of any $\Phi$,  the {\em{weight}} of $\Lambda$. It is independent of the 
choices.

\subsection{Pure modules.}  Fix an integer $j$. An indecomposable module 
$V_{\Lambda, t}$  for $K$ as above is {\em{$q$-pure of weight $j$}}, or simply 
{\em{pure}}, if \\

(i) the eigenvalues of $\Phi$ in $V_{\Lambda, t}$ are $q$-Weil numbers, and \\

(ii) $w(\Lambda)=t+j$.\\

\noindent By the argument at the end of the previous subsection, changing $\Phi$ will 
change the eigenvalues of $\Lambda(\Phi)$ only by  roots of unity, and hence both 
conditions are independent of the choice of $\Phi$.  Also, an indecomposable 
$V_{\Lambda,t}$ is $q_K$-pure of weight $j$ if and only if, for each finite extension 
$L$ of $K$, the restriction $V_{\Lambda,t}|_L$ of $V_{\Lambda,t}$ to $WD_L 
\subseteq WD_K$ is $q_L$-pure of (the same) weight $j$. To see this, note  since the 
condition is obviously stable under passage from $K$ to $L$, it is enough to show the 
descent statement from an $L$, as above, such that  $\Lambda|L$  is unramified. In this 
case, if $f=f(L/K)$ is the degree of the residue field extension, then $\Phi^f$ is a 
Frobenius element for $W_L$, and, in the above notation, $\chi(\Phi^f)= \chi_0(\Phi)^f$. 
Hence $\alpha =(\chi(\Phi^f))^{1/f}\zeta$, for some $f$-th root of $\chi(\Phi^f)$. 
Suppose now that $\chi(\Phi^f)$ is a  $q_L$-Weil number of weight j. Then, since 
$q_L=q_K^f$, $\alpha$ is a $q_K$-Weil number of weight j also. This shows (i) holds 
over $K$ if it holds over $L$.  To see (ii), just note that $w(\Lambda)$ is unchanged when 
$q_L=q_K^f$ is replaced by $q_K$ and $|\chi(\Phi^f)|$ is replaced by 
$|(\chi(\Phi^f))^{1/f}|$. This proves the claim. \\

\noindent We say that a general  member $V$ of $Rep^{s}(WD_K)$ is pure of weight 
$w$ ($w\in \mathbf{R}$) if each indecomposable constituent is pure of weight $w$. Of 
course, if the module $V$ is pure of weight $w$, then $w$ is  uniquely determined. 
Furthermore, if $V$ is pure of weight $w$, then any conjugate $\eta V$, for $\eta \in 
Aut(\mathbf{C})$ is also pure of weight $w$.\\

\noindent Finally, we say that an $l$-adic representation $V$ of $\Gamma_K$ is\\

\noindent (i) $q$-pure of weight $w$ if one, and hence any, associated member of 
$Rep^s(WD_K)$ is $q$-pure of weight w, and\\

\noindent (ii) pure if it is $q$-pure of weight w for some w. \\

\noindent Here are summarized some basic facts about elements of $Rep^s(WD_K)$ :\\

\noindent {\bf{Proposition 1}}

\noindent Let $V_1, ...., V_n $ be  $l$-adic representations of $\Gamma_K$. \\

\noindent (i)Let $V$ be the direct sum of the $V_i$. Then $V$ is $q$-pure of weight $j$ 
if and only if each $V_i$ is $q$-pure of weight $j$. \\

\noindent (ii) $V$ is $q$-pure of weight $j$ if and only if its contragredient $V^{*}$ is 
$q$-pure of weight $-j$.\\

\noindent(iii) if $V$ is $q$-pure of weight $j$, then the Tate twisted module $V(m)= V 
\otimes \chi_{l}^{m}$, where $\chi_{l}$  is the usual $l$-adic cyclotomic character, and  $m 
\in \mathbf{Z}$, is $q$-pure of weight $j-2m$.\\

\noindent(iv) If  $V$ and $W$ are $q$-pure of weights $k$ and $l$, their tensor product 
$V\otimes W$ is $q$-pure of weight $k+l$.  \\
\\

\subsection{Weight-Monodromy Conjecture.} This is the following 
statement(\cite{illusie}):

\begin {quote} Let $X$ be a projective smooth variety defined over the local field $K$. Let, 
as usual, for a rational prime $l$ which is different from the residue characteristic of $K$,  
$H^j(\overline{X}, \mathbf{Q}_l)$ be the $l$-adic \'{e}tale cohomology of 
$\overline{X}$, regarded as a $\Gamma_K$-module. Then the $\Gamma_K$-module  
$H^j(\overline{X}, \mathbf{Q}_l)$ is $q$-pure of weight $j$, where $q$ is the 
cardinality of the residue field of $K$. \end{quote}

\noindent Remark. Since $X$ is projective, its cohomology is polarizable, and so $\Phi$ 
acts on $det(H^j(\overline{X}, \mathbf{Q}_l))$  as $q^{jb_{j}/2}$
where $b_j$ is the dimension of $H^j(\overline{X}, \mathbf{Q}_l)$. On the other hand, if 
$H^j(\overline{X}, \mathbf{Q}_l)$  is $q$-pure of some weight $k$, we have also that 
$\Phi$ acts on this space as $q^{kb_{j}/2}$. Hence, we say simply that 
$H^j(\overline{X}, \mathbf{Q}_l)$ satisfies WMC if it is $q$-pure. \\

\noindent If, contrary to the convention of this paper, K and its residue field both have 
characteristic p, then WMC is a theorem of Deligne (\cite{weil2}, Theorem 1.8.4). In mixed 
characteristic, WMC is known for curves and abelian varieties (\cite{sga7-1}), for surfaces 
(\cite{rz}, Theorem 2.13, \cite{alter}, and see below, (4.3)), certain 
threefolds (\cite{ito3folds}) and, as mentioned in the introduction,  a class of  Shimura varieties associated to division algebras over CM fields  
(\cite{itopadic}, \cite{tywmc}).\\

\noindent  As  De Jong  remarks in the Introduction to \cite{alter}, it follows from his 
theory of alterations  that condition (i) of the definition of purity is always satisfied for l-
adic representations that are subquotients of the $l$-adic \'etale cohomology of a 
quasiprojective-variety $X$ over a non-archimedian local field $K$.  We sketch this result, 
for the case that $X$ is smooth and projective, since it is basic.\\

\noindent {\bf{Proposition 2.}}(De Jong). Let $X$ be a smooth projective variety defined 
over the local field $K$. Let $\Phi \in \Gamma_K$ be a Frobenius. Then the eigenvalues of 
$\Phi$ on  $H^j(\overline{X}, \overline {\mathbf{Q}_l})$ are integral Weil numbers. \\

\noindent {\bf{Proof.}}  Let $L$ be a finite extension of $K$ over which there exists an 
$L$-alteration $a: X^{\prime} \rightarrow X_L$ such that $X^{\prime}$ is the generic 
fiber of a strictly semistable scheme $\cal{X}^{\prime}$ defined  over the ring of integers 
$\cal{O}_L$ of $L$. Since an alteration is surjective and generically finite,  we may regard  
$H^j(\overline{X}, \overline {\mathbf{Q}_l})$ as a submodule of 
$H^j(\overline{X^{\prime}}, \overline {\mathbf{Q}_l})$ via $a^*$. Let 
$\underline{\cal{X}^{\prime}}$ be the geometric special fiber of $\cal{X}^{\prime}$. Since  
$\cal{X}^{\prime} $ is strictly semistable, its cohomology is computable via the 
$\Gamma_L$-equivariant weight spectral sequence of Rapoport and Zink (c.f.\cite{rz}, 
Section 2, and \cite{illusie},3.8). In the notation of \cite{illusie}, we have

$$ _{W}E_1^{ij}=H^{i+j}(\underline{\cal{X}^{\prime}}, gr^W_{-
i}R\Psi(\overline{\mathbf{Q}_l})\Longrightarrow H^{i+j}(\overline{X^{\prime}}, 
\overline{\mathbf{Q}_l}). $$

\noindent Thus, it suffices to show that the eigenvalues of $\Phi$ on each  
$_{W}E_1^{ij}$ are integral Weil numbers. But each  $_{W}E_1^{i,j}$ is a direct sum 
of cohomology groups of the form

$$H^{j+i-2l}(\underline{{\cal{X}^{\prime}}}^{(2l-i+1)}, \overline{\mathbf{Q}_l})(i-
l)$$ 
\noindent where $l \geq max (0, i)$, and  $\underline{\cal{X}^{\prime}}^{(2l-i+1)}$ is the disjoint union of smooth 
proper subvarieties of  $\underline{\cal{X}^{\prime}}$ defined by taking   $(2l-i+1)$-
fold intersections of  the irreducible divisors provided in the definition of the strict 
semistability of $\cal{X}^{\prime} $. The result now follows from the Weil conjectures.

\subsection{Weight-Monodromy: Background Facts}

Let $W$ be a finite set of $q$-Weil numbers and $m: W \rightarrow \mathbf{Z}_{\geq 
0}$ be a function. For us,  $m(\alpha)$ is the multiplicity of $\alpha$ in the spectrum of a 
Frobenius in a semistable module.  The pair $(W, m)$ is said to be {\em{wm-$q$-pure of 
weight j}} if , \\

\noindent (i) whenever $|\alpha| > q^{j/2}$, $m(q^{-1}\alpha)\geq m(\alpha)$, and,\\

\noindent (ii) for all $\alpha$, $m(\alpha)=m(q^{-s_{\alpha}}\alpha)$, where 
$s_{\alpha}=2log_{q}(|\alpha|q^{-j/2})$.\\

\noindent Let, for $\alpha \in W$,  $|\alpha| \geq q^{j/2}$, $$\delta(\alpha)=m(\alpha)-
m(q\alpha).$$

\noindent Let $W^{+}$ be the subset of $W$ of all $\alpha$ such that $|\alpha| \geq 
q^{j/2}$.\\

\noindent If  $(V, \rho)$ is a $q$-pure of weight j semistable representation of $WD_K$, let 
$b(V)$ denote the number of indecomposable factors in any representation of $V$ as a 
direct sum of such so that $b(V)=dim(ker(N_{\rho}))$. More generally, for any nilpotent 
endomorphism $N$ of $V$, let $b(N)=dim(ker(N))$ denote the number of 
indecomposable Jordan blocks in the representation of $N$ as a direct sum of such. 
Evidently,  if  $V$ is a $q$-pure of weight j semistable representation of $WD_K$, the 
associated pair $(W_{V}, m_{V})$ is wm-$q$-pure of weight j.  Conversely, if $(W, m)$ is 
wm-$q$-pure of weight j, then

$$ \bigoplus_{ \alpha \in W^{+}}V_{\alpha, s_{\alpha}}^{\delta(\alpha)}$$

\noindent  belongs to the   unique isomorphism class of $q$-pure of weight $j$ semistable 
representations $(V_{(W,m)}, \rho_{(W,m)})$ of $WD_K$ that give rise to $(W, m)$.

\noindent In this case let $$b(W, m)=b(V_{W})=\Sigma_{\alpha \in 
W^+}\delta(\alpha)$$.

\noindent The following elementary result is key to our work in this paper.\\

\noindent {\bf{Proposition 3.}} Let $K$ be a local field and let $V$ be a finite 
dimensional representation in $Rep^s(WD_K)$. Let $F^{\cdot}V$ be a filtration of 
$V$ by $WD_K$-stable submodules.  Suppose that the graded Galois module 
$Gr_{F}(V)$ is $q$-pure of weight $j$. Then $V$ is $q$-pure of weight $j$.

\noindent {\bf{Proof}.} Restricting from $K$ to a suitable extension $L$, we can assume 
that $V$ is semistable. \\

\noindent Let $E_{V}$ be an endomorphism of a vector space $V$, and suppose that we 
have an $E_V$ stable short exact sequence

$$ 0\rightarrow S \rightarrow V \rightarrow Q\rightarrow 0$$

\noindent  with  induced endomorphisms $E_S$ and $E_Q$ on $S$ and $Q$. Let 
$K_S$, $K_V$, and $K_Q$ be the kernels of these operators. Then

$$dim(K_{S}) + dim(K_{Q}) \geq dim(K_{V}).$$

\noindent  This is evident since  we have a short exact sequence

 $$ 0\rightarrow K_{V} \cap S \rightarrow K_{V} \rightarrow \frac{K_{V} + 
S}{S}\rightarrow 0$$

\noindent and $K_{V} \cap S=K_{S}$ and $\frac{K_{V} + S}{S}$ is a subspace of 
$K_{Q}$.

\noindent Hence, by induction,  if $E_V$ is a filtered endomorphism of  $F^{\cdot}V$, 
inducing $Gr_{F}(E_{V})$ on  $Gr_{F}(V)$ , then

$$dim(ker(Gr_{F}(E_{V}))) \geq dim (ker(E_{V})).$$

\noindent  We apply this to the case $E_V=N$.\\

\noindent{\bf Lemma}. Let  $(W, m)$ be wm-$q$-pure of weight j. Let $(V, \rho)$ be a 
semistable representation of $W_K$ such that $(W_{V}, m_{V})=(W, m)$. Then

$$ b(V) \geq   b(W, m). $$

\noindent Further,

$$ b(V) =b(W, m)$$

\noindent if and only if $(V, \rho)$ is $q$-pure of weight $j$.\\

\noindent {\bf {Proof.}} Obvious. \\

\noindent To conclude the proof of the Proposition, we note that $Gr_{F}(V)$  defines the 
same pair $(W, m)$ as $V$. Since we assume $Gr_{F}(V)$ is pure, we have $b(W, 
m)=b(Gr_{F}(V))$. By the remarks just above,  we always have  $b(Gr_{F}(V)) \geq 
b(V)$ and $b(V) \geq b(W,m)$. Hence $b(V)=b(W,m)$ and so $V$ is $q$-pure of 
weight j.  \\

\subsection{ A problem on abelian varieties.}

\noindent {\bf{Proposition 4.}}

\noindent Let $A$ be an abelian variety defined over a number field $J$. Let $M$ be an 
irreducible motive, defined over $J$ in the category of motives for absolute Hodge cycles 
generated by $A$(\cite{tancats}). Then for each prime $l$ and each finite place $v$ of 
$J$,  the $l$-adic cohomology $M_l$ of $M$ satisfies the WMC.\\

\noindent {\bf{Proof}}. This is, of course, trivial: any irreducible $M$ is of then form 
$M_0(n)$ where $M_0$ is a submotive of the motive $\otimes^{k} H^1(A, 
\overline{Q})$, $k$ is a non-negative integer,  and $n$  is the $n$-fold Tate twist. The $l$-adic cohomology $M_{0,l}$ of $M_0$ is a  
$Gal(\overline{J}/J)$ direct summand of $\otimes^{k} H^1(A, \overline{Q}_l)$ and 
hence everywhere locally satisfies WMC since $H^1(A, \overline{Q}_l)$ does. \\

\noindent Problem: Is the conjugacy class of $N_l$(in $GL(M_l)$)  independent of $l$? 
Evidently, this amounts to asking whether the  Frobenius eigenvalues on the 
semisimplification of $M_l$ is independent of $l$. Of course, these statements are 
consequences of the standard l-independence conjecture of Serre and Tate which asserts 
that, for any motive $M$ over $J$, and any non-archimedian completion $J_v =K$, the 
isomorphism classes of the elements of $Rep^{s}(WD_K)$ gotten from the $l$-adic 
\'{e}tale cohomology groups of $M$ are all the same.\\

\section{Automorphic Forms.}

\subsection{Basic Conventions}

\noindent Let $F$ be a number field with adele ring $\mathbf{A}_F$. Let 
${\emph{A}}_{0}(F, n)$ be the set of irreducible cuspidal unitary summands of the space 
$L^2(GL(n,\mathbf{A}_F)/GL(n,F))$. Each constituent $\pi$ of such a space is 
isomorphic to a restricted tensor product $\pi = \otimes_v \pi_v$ where $\pi_v$ is an 
infinite dimensional irreducible unitary representation. If $v$ is finite, each $\pi_v$ is 
classified up to isomorphism by an associated isomorphism class $\sigma(\pi_v)$ of n-
dimensional members of  $Rep^{s}(WD_v)$, where we denote by $WD_v$  the Weil-
Deligne group of $F_v$ (\cite{ht},\cite{kull}). As is customary, we denote also by  
$\sigma(\pi_v)$ any member of its class. Let $W_v$ be the Weil group of $F_v$. Then 
$\sigma(\pi_v)$ is isomorphic, as in  1.9, to a direct sum of indecomposable modules of the form
$V_{\Lambda_i, t}=\Lambda_i \otimes V_{q^{t/2}, t}$,  $1 \leq i \leq n_v$, with 
irreducible representations $\Lambda_{i}$ of $W_v$. \\

\subsection{Ramanujan Conjecture.} This is the assertion:

\begin{quote}

Let  $\pi \in \emph{A}_{0}(F, n)$. Let $v$ be a finite place of $F$ and  define, as above,  
the set of representations $\Lambda_i$ of $W_v$ for  $\pi_v$. Then the  image of each 
$\Lambda_i$ is bounded. \end{quote}

\noindent {\bf Remark}: This form of the conjecture is equivalent to the more elementary 
statement, independent of the Local Langlands Correspondence,  which asserts that each 
$\pi_v$  is tempered.  However, we work exclusively with the formulation via Weil-Deligne 
groups in this paper.\\

\noindent Suppose now that $n=2$. Then, at non-archimedian $v$,  the local components 
$\pi_v$ of a cuspidal $\pi$ are classified into several types: \\

\noindent (i) $\pi_v$  is supercuspidal, \\

\noindent (ii) $\pi_v$  is a twist of the Steinberg representation: $\pi_v= St_v\otimes 
\psi(det)$, so that $\sigma(\pi_v)=\psi \otimes V_{q^{1/2},1}$ with a character  $\psi$ of 
$F_v^{*}=W_v^{ab}$.\\

\noindent (iii) $\pi_v$  is  principal series.\\

\noindent In cases (i) and (ii), $\sigma(\pi_v)$ is indecomposable and if $\pi$ is unitary, $\Lambda=\Lambda_1$  is bounded. 
(In case (i), $\Lambda$ is irreducible, $t=0$, and $det(\Lambda)$ is the unitary central 
character of $\pi$, so $\Lambda$ is bounded; in case (ii), $t=1$, so $\Lambda =\psi$ is 
one-dimensional, and $\Lambda^2$ is the unitary central character of $\pi$, so $\Lambda$ 
has bounded image. \\

\noindent For case (iii), $\sigma(\pi_v)$  is a direct sum of 2 quasicharacters $\psi_1$ and 
$\psi_2$ of $F_{v}^*$ whose product is the central character of $\pi$, hence unitary. The 
classification of unitary representations shows that  either (a) $|\psi_1| = |\psi_2|=1$ or (b) 
there are quasicharacters $\mu=|\cdot|^t$ where $0<t<1/2$ and $\psi$ of $F_{v}^{\ast}$ 
such that $\sigma(\pi_v)$ is the sum of quasicharacters $\mu\psi$ and $\mu^{-1}\psi$.  
Hence, the Ramanujan Conjecture amounts to the assertion that  representations of this type 
({\em{complementary series}}) don't occur as local components of cusp forms. Note that 
at such a place, the local central character $\omega_{\pi, v}$ of $\pi$ is $\psi^2$. For the 
forms of interest in this paper,  $F$ is totally real, and  the infinity type $\pi_{\infty}$ of 
$\pi$ is discrete series and has the property that the idele class character $\omega_{\pi}$ 
takes the form $\omega_{\pi}=\nu_F^j \otimes \phi$, where $\nu_F $ is the norm, $j$ is an 
integer, and $\phi$ is a character of finite order.  Hence, $\psi^2=\phi_v$ and so $\psi$ 
has finite order.  Invoking the Gruenwald-Hasse-Wang theorem, we see that there is an 
idele class character of finite order $\eta$ such that the local identity $\eta_v= \psi$ holds. 
Thus, replacing $\pi$ by a form of the same infinity type $\pi^{\prime}=\pi \otimes \eta^{-
1}$, we see that to establish the conjecture for all local components of all cusp forms 
$\pi^{\prime}$ of the given discrete series infinity type, it is enough to prove it for all $\pi$ 
of the given type at all $v$ that are unramified for $\pi$. Although this easy argument is 
special to $GL(2)$, it may be worth noting that solvable base change for $GL(n)$ should 
provide a reduction of Ramanujan to the case of semistable representations (i.e. to those 
whose local components $\sigma(\pi_v)$ are semistable.) 
 
\subsection{Global Langlands Correspondence.}

Let  $F$ be a number field and let $(V_{l}, \rho_l)$ be an irreducible n-dimensional $l$-adic 
representation of $\Gamma_F$. Fix, for the rest of the paper, an isomorphism $\iota_{l}: 
\overline{\mathbf{Q}_l}\rightarrow  \mathbf{C}$. For each finite place $v$ of $F$, 
whose residue characteristic is different from $l$,  choose $t_l$ and $\Phi_l$, as before. 
Let $\rho_{l, v}^* $ be the associated member of $Rep^{s}(WD_v)$  so defined. \\

\noindent {\bf{Global Langlands Correspondence(GLC).}} This is the assertion:

\begin{quote}

Suppose that the irreducible $l$-adic representation $(V_{l}, \rho_l)$ is motivic. Then there 
are cuspidal representations  $\pi \in {\emph{A}}_{0}(F, n)$ and $\chi \in 
{\emph{A}}_{0}(F, 1)$   such that,  for all $v$ whose residue characteristic is different 
from $l$,  $\sigma(\pi_v)\otimes \chi_v $ is the class of  $\rho_{l, v}^*$ . \end{quote}

\noindent Remark 1. Since the statement of the GLC presupposes the existence of the Local 
Langlands Correspondence and an $l$-adic representation of the absolute Galois group of a global field,  the GLC is often called the problem of {\emph {Local-Global Compatibility}.\\

\noindent Remark 2. If the residue characteristic of $v$ is $l$ the  classes $\sigma(\pi_v)\otimes \chi_v $  
can be predicted using methods of $p$-adic cohomology ([Fo]). This done, the above 
conjecture is extended to all finite places.\\

\noindent Remark 3. It is usual, especially to treat compatibility questions as l and $\iota_l$ 
vary, to formulate the conjecture in terms of a motive M and its Galois representations. 
However, as we do not treat compatibility questions in any essential way in this paper, there 
is no benefit to this viewpoint. \\

\noindent Remark 4. There is a converse conjecture: if  $\pi_{\infty}$ is algebraic ([C1]) 
then there should exist a $(V_{l}, \rho_l)$ corresponding to $\pi$ as above.

\subsection{GLC and WMC}

\noindent {\bf{Proposition 5}}

\noindent Suppose that the GLC holds for the motivic $l$-adic representation $(V_{l}, 
\rho_l)$ over $F$. Then WMC holds for $(V_{l}, \rho_l)$.\\

\noindent {\bf{Proof.}}

\noindent The conjecture is invariant under Tate twist, so we may assume  that $(V_{l}, 
\rho_l)$ is isomorphic to a subquotient of $H^{i}(X, \overline{\mathbf{Q}}_l)$ for some 
smooth projective $X$ over $F$. For almost all places $v$, $\pi_v$ is  unramified. At such 
a place, the parameter $\sigma(\pi_v)$  consists of $n=dim(V_l)$  unramified 
quasicharacters of $F_{v}^*$, whose values on a prime element of $F_v$ determine the 
unordered n-tuple $\{ \alpha_j | j=1, ...,n\}$. Each $\alpha_j$ is a Weil number, and, if we 
further restrict $v$ to be a place of good reduction of $X$, then we have 
$|\alpha_j|=q^{i/2}$ for all $j$, for some $i$ which is independent of $j$. Consider the 
cuspidal representation $\pi^{\prime}= \pi\otimes |\cdot|^{i/2}$. Then $\pi^{\prime}$ is 
unitary because its central character is unitary; this holds at all unramified places $v$ and 
hence everywhere.  Let $v_0$ be finite place which we wish to study. The classification 
(\cite{tadic}, see \cite{kull})of unitary representations of $GL(n, F_{v_0})$ shows that  
$\sigma(\pi^{\prime}_{v_0})$  is a direct sum of  indecomposables  $\Lambda \otimes 
V_{q^{t/2}, t}$ where $$ -1/2 < w(\Lambda) < 1/2.$$

\noindent Hence

$$ (i-1)/2 < w(\Lambda \otimes |\cdot|^{-i/2})<(i+1)/2.$$

\noindent  Since $(V_{l}, \rho_l)$ is motivic, Proposition 1 shows that $w(\Lambda 
\otimes |\cdot|^{-i/2})$ is an integer in this interval. Hence $w(\Lambda \otimes |\cdot|^{-
i/2})=i$, which means $w(\Lambda)=0$, as was to be shown.\\

\noindent Remark. The proof of Proposition 5 uses only the fact that the weight of 
$\Lambda$ is an integer, not the fact that the eigenvalues  of Frobenius are algebraic 
numbers.

\section{Zeta functions of quaternionic Shimura varieties.}

\noindent Assume henceforth that $F$ is totally real and let $G$ be an inner form of 
$GL(2)/F$, so that $G(F)=B^*$ with a quaternion algebra $B$ over $F$.  Let $J_{F, 
nc}= \{\tau_{1} , ..., \tau_{r}\}$ be the set of real embeddings (= infinite places) of $F$ 
where $B$ is indefinite; assume that  $J_{F, nc}$ is non-empty and contains 
$\tau_{1}=1_F $. To $B$ is attached a Shimura variety  $Sh_B$ defined over 
$F^{\prime}$, the smallest extension of $\mathbf{Q}$ containing all  elements 
$\tau_{1}(f) +...+ \tau_{r}(f)$ for all $f \in F$. See (\cite{shcurves}, \cite{dsv1}) for 
constructions of  $Sh_B$. It is the projective limit of quasi-projective r-folds  $Sh_{B, 
W}$, where $W$ is an open compact subgroup of  $G(\mathbf{A}_{F,f})=(B\otimes \mathbf{A}_{F,f})^{*}$. Each  
$Sh_{B, W}$ is defined over $F^{\prime}$, and is a finite disjoint union of connected r-
folds. These components are proper if $G$ is not $GL(2)/F$ and smooth if $W$ is small 
enough.  Any such Shimura variety is called a {\em{quaternionic Shimura variety}}.  The  
Hasse-Weil zeta function, at almost all places, of the $l$-adic \'etale cohomology of such 
Shimura varieties has been computed by Reimann (See \cite{ssz}, Theorem 11.6), in the 
case $B\neq GL(2)/F$,  and Brylinski and Labesse in the case $G=GL(2)/F$.  In the latter 
case, it is the zeta function of an intersection cohomology which has been computed, and it 
is this cohomology that we consider in the following, using the same notation as the other 
cases. The zeta functions of the $l$-adic cohomology groups of $Sh_{B, W}$ have, at almost 
all places of $F^{\prime}$, the form conjectured by Langlands (\cite{zssv})and proved by 
him in the case $r=[F:\mathbf{Q}]$ (\cite{zssv}). See (\cite{zeta}, Sections 3.5, 5.1, and 
7.2) for an expository treatment of the result but not the proof.\\

\noindent  For our purposes, it is sufficient to give a global description of the result over a 
Galois extension $L$ of $\mathbf{Q}$ which contains $F$. Thus, for each $j \in \{1,..., r 
\}$, let $ \overline{\tau_j}$ be an extension of  $\tau_j$ to $\overline{L}$. Let $\pi$  be a cuspidal 
automorphic representation of weight $(2, ..., 2)$ of  $GL(2, \mathbf{A}_{F})$ which is 
discrete series at any finite place of $F$ at which $B$ is ramified. Choose $\pi$ so that its 
central character $\omega_{\pi}$ satisfies $\omega_{\pi}=\Psi |\cdot|^{-1}$ with a 
character $\Psi$ of finite order.  Let $T$ be the number field generated by the traces 
$tr(\sigma(\pi_v)(\Phi))$ for all $v$ which are unramified for $\pi$. As shown by Taylor 
(\cite{rthmf}, \cite{rthmf2}), there is an irreducible two-dimensional $l$-adic representation $\rho_l^T$, 
depending only on  $\pi$ and $\iota_l$,  which satisfies GLC relative to $\pi$. Let 
$\rho^T_{l, L}$ be the restriction to $\Gamma_L$ of $\rho^T_l$. Let 
$^{[\overline{\tau_{j}}]}\rho^T_{l, L}$  be the representation defined, for $\eta \in 
\Gamma_L$, by

$$ ^{[\overline{\tau_j}]} \rho^T_{l, L}(\eta)=\rho^T_{l, L}({\overline{\tau_j}}\eta 
{\overline{\tau_j}}^{-1}).$$

\noindent Let $$R_l(\pi)= R_{l,J_{F, nc}}(\pi)= \otimes_{j=1,...,r} {^{[\overline{\tau_j}]} 
\rho^T_{l, L}}.$$

\noindent Then $R_l(\pi)$ is a semisimple $l$-adic representation of $\Gamma_L$ of 
dimension $2^r$.

\subsection{Semisimple cohomology of $Sh_B$.} Let, if it exists,  $\pi^{\prime}$ be an 
automorphic representation of $G(\mathbf{A}_{F})$ such that $\pi^{\prime}_v$ is 
isomorphic to $\pi_v$ at all places $v$ of $F$ which are unramified for $B$. Thus $\pi= 
JL(\pi^{\prime})$, where JL denotes the Jacquet-Langlands correspondence. Choose an 
open compact subgroup $W$ as above so that $Sh_{B, W}$ is smooth and  
$\pi^{\prime}_f $ has a non-zero space of $W$-invariants. Let $\pi{^{\prime}}_{f,W} $ 
denote the representation of the level $W$ Hecke algebra ${\cal{H}}_{W}$ gotten from 
$\pi^{\prime}_f $. Let

 $$H^r(Sh_{B, W}, \overline{\mathbf{Q}_l})(\pi{^{\prime}}_{f,W})$$

\noindent be the $\pi{^{\prime}}_{f,W}$-isotypic component of $H^r(Sh_{B, W}, 
\overline{\mathbf{Q}_l})$.\\

\noindent {\bf{Proposition 6}}

The irreducible subquotients of  the action of $\Gamma_L$ on $H^r(Sh_{B, W}, 
\overline{\mathbf{Q}_l})(\pi{^{\prime}}_{f,W})|L$

are exactly the irreducible subquotients of  $R_l(JL(\pi^{\prime}))$.\\

\noindent {\bf Proof.} By the $l$-adic Cebotarev Theorem (\cite{abladic}), it suffices to 
show that the semisimplification of the Galois action on $H^r(Sh_{B, W}, 
\overline{\mathbf{Q}_l})(\pi{^{\prime}}_{f,W})|L$ is a multiple of 
$R_l(JL(\pi^{\prime}))$. But, up to notation and the base change to $L$, this is given by 
Theorem 11.6 of \cite{ssz}, and by the main theorem of \cite{brylab} 
in the non-compact case. See Section 5.3 for some explicit review of the zeta function.

\section{Ramanujan and Weight-Monodromy for Hilbert modular forms}

\noindent Let  $F$ be a totally real field and let $\pi=\pi_\infty\otimes \pi_f$  be a 
holomorphic cuspidal automorphic representation of $GL(2, \mathbf{A}_F)$. Up to twist, the 
isomorphism class of $\pi$ at the infinite places of $F$ is specified, as usual,  by a tuple of 
positive integral weights $k=(k(\tau))$, where the variable $\tau$ runs over the real 
embeddings of $F$. We normalize $\pi$ by assuming that its central character 
$\omega_{\pi}$ satisfies $|\omega_{\pi}|= |\cdot|^{1-k}$ where $k$ is the maximum of 
the $k(\tau)$'s. It is natural to  classify the holomorphic cuspidal $\pi$'s into several types, 
depending on $\pi_{\infty}$, i.e. on the classical weights at each infinite place: \\

\noindent (i) type G: all the weights are 1. \\

\noindent (ii) type MC: all the weights are at least 2 and they are all congruent modulo 2;\\

\noindent (iii) type NMC: all the weights are at least 2 and they are not all congruent 
modulo 2;\\

\noindent (iv) type NC: at least one, but not all , of the weights are 1.\\

\noindent Types G and MC are well-studied.  RC is known at all places for type G 
(\cite{ordinary}, \cite{rogtun}); there is a 2-dimensional Artin representation   $\rho$ of 
$\Gamma_F$ that satisfies the GLC. In this paper we prove RC at all places for the class of 
forms MC. As mentioned in the Introduction, the method of this paper should apply to type 
NMC but we do not consider this case in the paper. Type NC, except for the case of CM 
forms of this type, is completely open. Even in the case where the weights are  all congruent 
mod 2, we do not know any motivic realization of associated Galois representations 
(\cite{jarvis}).

\subsection {Proof of Theorem 1.} \noindent We must show that the  Ramanujan 
Conjecture (RC) holds for all Hilbert modular forms of type MC. \\

\noindent   Let $\pi$ be a representation of type MC of classical weight $k=(k(\tau))$. Let 
$T$ be the field generated by almost all Hecke eigenvalues of $\pi$. Let $\rho_l^T$ be one 
of the  $[T:\mathbf{Q}]$ two dimensional  $l$-adic representations attached to $\pi$ which 
satisfy GLC.  As shown in \cite{brhmf},  these representations are motivic except possibly 
in the case where $[F:\mathbf{Q}]$ is even, $k(\tau)=2$ for all $\tau$, and $\pi_v$ 
belongs to the principal series for all finite $v$.  Hence, except in this case, RC follows 
from Proposition 5.\\

\noindent Let $v$ be a finite place at which we will prove that $\pi_{v}$ satisfies 
RC. Changing $l$, if necessary, we assume that $v$ does not lie above $l$. Replacing $\pi$ by a twist $\pi\otimes \Psi$, we may assume that $\pi$ is unramified 
at all finite places of $F$ which lie above the rational prime $p$ under $v$ . Let 
$\tau_1$ be the tautological infinite place of $F$ and  let $\tau_2$ be another infinite place. 
Let $B$ be the quaternion algebra over $F$ which is unramified at $\tau_1$, $\tau_2$ , 
and at all finite places, and which  is ramified at the remaining infinite places. Let $G$ be 
the inner form  of $GL(2)$ over $F$ such that $G(F)=B^{\ast}$. Let $L$ be a Galois 
extension of $\mathbf{Q}$ which contains $F$. By Proposition 6, the 4-dimensional l-
adic representation $R_l(\pi)$ of $\Gamma_L$, made using $\rho^{T}_{l}|_{L}$ with 
$J_{F,nc}= \{\tau_1, \tau_2 \}$, is isomorphic to the sum of  irreducible subquotients of 
the cohomology $H^2(Sh_{B, W}, \overline{\mathbf{Q}_l})$ of the quaterionic Shimura surface 
$Sh_{B, W}$, for small enough open compact subgroup $W$ of the finite adele group 
$(B\otimes \mathbf{A}_{F,f})^{*}$.\\

\noindent Now we need to make explicit the action on $R_l(\pi)$ of a decomposition group 
of a place $w$  of $L$ dividing $v$. Choose a decomposition group 
$D_{\overline{w}} \subset \Gamma_F$ for $w$, and denote by $R_{l, 
w}(\pi)$ the restriction of $\rho^{T}_{l}|_{L}$ to $D_{\overline{w}}\cap \Gamma_L$. Let $\tau_{2} v$ be the place of $F$ lying below $\overline{\tau_2}w$. Let $f_1$ and 
$f_2$ be degrees of the residue field extensions associated to $L_w / F_v$, and 
$L_{\overline{\tau_2}w}/F_{\tau_{2}v}$, respectively. For each place $v^{\prime}$ of 
$F$ above $p$, the  restriction of $\rho^T_l|L$ to $D_{v^{\prime}}$ is unramified. Denote 
the eigenvalues of  $\Phi_v$ by  $\alpha_v$ and $\beta_v$. Denote the eigenvalues of  
$\Phi_{\tau_{2}v}$ by  $\gamma_{\tau_{2}v}$ and $\delta_{\tau_{2}v}$.  Over $L$, 
$\Phi_w$ acts via $\rho^T_{l}| L$  with eigenvalues $\alpha_v^{f_1}=\alpha_w$ and 
$\beta_v^{f_1}=\beta_w$.  Likewise,  $\Phi_{\overline{\tau_2}w}$ acts with eigenvalues 
$\gamma_{\tau_{2}v}^{f_2}= \gamma_{\overline{\tau_2}w}$  and $\delta_{{\tau_2} 
w}^{f_2}=\delta_{\overline{\tau_2}w}$. Hence $\Phi_w$ acts via  
$^{[\overline{\tau_2}]}\rho^{T}_{l}| L$ with eigenvalues  
$\gamma_{{\overline{\tau_2}}w}$  and  $\delta_{{\overline{\tau_2}}w}$. Note  the 
product relation $\gamma_{{\overline{\tau_2}}w}\delta_{{\overline{\tau_2}}w}=\zeta 
q_{\overline{\tau_2}w}= \zeta q_w$, where $\zeta$ is a root of unity,  which is obvious 
since, by the global Langlands correspondence,  $det(\rho^T_l(\Phi_v))= \mu q_v$ with a 
root of unity $\mu$, where $q_v$, $q_w$, and $q_{\overline{\tau}_{2}w}$ are the numbers of elements of the residue fields associated to $v$, $w$, and $\overline{\tau}_{2}w$.  \\

\noindent By definition of  $R_l(\pi)$, $\Phi_w$ acts via $R_l(\pi)$ with eigenvalues 
$a=\alpha_w \gamma_{\overline{\tau_2}w}$,  $b=\alpha_w \delta_{\overline{\tau_2}w}$, 
$c=\beta_w \gamma_{\overline{\tau_2}w}$,  and $d=\beta_w 
\delta_{\overline{\tau_2}w}$.\\

\noindent By Proposition 2, $a, b, c$ and  $d$ are $q_w$-Weil numbers. Hence 
$ab=\alpha_{w} ^2 \zeta q_w$, and so $\alpha_{w} ^2 $ is a $q_w$-Weil number. Since 
$\alpha_{w} ^2 =(\alpha_{v} ^2)^{f_1}$, we see that  $\alpha_{v} ^2$ is a $q_v$-Weil 
number. Likewise, $\beta_{v} ^2$ is a $q_v$-Weil number. If we let $|\alpha_{v} 
^{2}|=q_{v}^{l/2}$ and $|\beta_{v} ^{2}|=q_{v}^{m/2}$ then 
$|\alpha_{v}|=q_{v}^{l/4}$ and $|\beta_{v}|=q_{v}^{m/4}$, for integers $l$ and $m$. \\

\noindent Recall now the following (\cite{sharam}) Ramanujan estimate:

$$  q_{v}^{-1/5}< |\alpha_{v}|, |\beta_{v}| < q_{v}^{1/5},$$

\noindent which applies to all unitary cusp forms $\pi^{\prime}$ for $GL(2)$ over any 
number field, at an unramified place $v$ for $\pi^{\prime}$. Since in our case $\pi \otimes 
|\cdot|^{1/2}$ is unitary, we see that the exponents $q_{v}^{\frac{l-2}{4}}$ and 
$q_{v}^{\frac{m-2}{4}}$ must be compatible with this estimate. Evidently this happens if 
and only if $l=m=2$, which is precisely what we needed.\\

\noindent Of course, this result about $\pi$  implies something about $\rho^{T}_{l}$: \\

\noindent {\bf{Corollary 7.}} Let $\pi$ be a Hilbert modular form of type MC. Then any l-
adic representation $\rho^T_l$ which satisfies the GLC satsifies WMC at all places $v$ 
prime to l. \\

\noindent Remark.  An easy extension of the above method shows RC at all places for all $F$, at least under the congruence condition on the weights. To prove RC at the place $v$, it is enough to choose any totally real quadratic extension $K$ of $F$. Then, defining B over $K$ as above, one proves  RC, by the method here,  at each place of $K$ for the base change $\pi_K$ of $\pi$ from $F$ to $K$. But it is easy to see that RC holds for $\pi_K$ at a place $w$ of  $K$ iff it holds for  $\pi$ at the place of $F$ under $w$. Thus,  RC may be proved for all Hilbert modular forms which satisfy the congruence condition  at infinity by a uniform method which reduces the problem to the calculation of \cite{ssz} and Shahidi's estimate. \\

\noindent {\bf{Proof of the Corollary}}. Indeed, it only remained, in view of the work of Carayol 
(\cite{hmfzeta}), to establish the result at the unramified places of $\pi$, and this is 
precisely the RC.

\subsection{Geometric Proof of Theorem 1.} \noindent 
 There exists a finite extension $\overline{L}_{\overline{u}}$ of $L_w$ over which the generic fiber $\overline{Sh_{B, W}}$ of a semistable alteration  of $Sh_{B, W}$ is defined. Then  $H^2(Sh_{B, W}, \overline{\mathbf{Q}_l})$ is direct summand, as $D_{\overline{u}}$-module of  $H^2(\overline{Sh_{B, W}}, \overline{\mathbf{Q}_l})$. This latter group satisfies WMC by \cite{rz}. Now $R_l(\pi)$ is, after some base change, a tensor product, and its associated Weil-Deligne parameter is thus a tensor product as well.  
We now note the  following simple result whose proof is left to the reader: \\

 \noindent {\bf{Lemma}}.  Let $V_1$ and $V_2$ be 2-dimensional representations of a Weil-Deligne group $WD$. Suppose that $V_1\otimes V_2$ satisfies WMC, and suppose that the modules $\Lambda^2(V_i)$ are pure of weight 2. Then each $V_i$ is pure of weight 1. \\

 \noindent Applying this with  $V_1 = \rho_l^T$, we conclude that  $\rho_l^T|_{D_{\overline{u}}}$ is pure of weight 1 at $u$, and hence $\rho_l^T$ is pure of weight 1 at $v$. Hence, since Local-Global Compatibility is known (\cite{rthmf}),  $\pi_v$ satisfies RC.

\section{Weight-Monodromy Conjecture for Quaternionic Shimura Varieties.}

\subsection{Proposition 8.} Let F be a number field and let V be a variety defined over 
$F$. Let l be a rational prime with $(v, l)=1$. Let $v$ be a finite place of $F$. Then if the 
WMC holds at the place $v$ for the semisimplification of $H^j(X, 
\overline{\mathbf{Q}_l})$ as a $\Gamma_F$-module, then the WMC holds for 
$H^j(X,\overline{\mathbf{Q}_l})$ at $v$.\\

\noindent {\bf{Proof.}} This is just a geometric restatement of a special case of Proposition 
3.

\subsection {Proof of Theorem 2:WMC for Quaternionic Shimura varieties.}

\noindent Let, with notations as above, $Sh_{B, W}$ be a quaternionic Shimura variety of 
dimension r. The Hecke algebra at level W acts semisimply, and  $H^{r}(Sh_{B, W}, 
\overline{\mathbf{Q}_l})$ is thus a direct sum of isotypic components $H^{r}(Sh_{B, 
W}, \overline{\mathbf{Q}_l})(\pi^{\prime}_{f, W})$. It is sufficient, in view of 
Proposition 1, to show that WMC holds for each of these components.  If $\pi_f $ is one-dimensional,   the result of Reimann (\cite{ssz}, Theorem 11.6) shows that after a finite base change to a number field $L$, the character of the Galois action on $H^{*}(Sh_{B, W}, \overline{\mathbf{Q}_l})(\pi_{f, W})$ is a sum of powers of the cyclotomic character at almost all, and hence, by Cebotarev,  all, finite places. Since $H^{j}(Sh_{B, W}, \overline{\mathbf{Q}_l})(\pi_{f, W})$ is pure of weight $j$ at almost all finite places $v$, the Galois action on it over $L$ is a multiple of $ \chi_l^{-j}$ where $\chi_l$ is the $l$-adic cyclotomic character.  Hence $H^{j}(Sh_{B, W}, \overline{\mathbf{Q}_l})(\pi_{f, W})$ is pure of weight $j$ at all places. (We note that this fact is much less deep: it is not hard to see that all of $H^{j}(Sh_{B, W}, \overline{\mathbf{Q}_l})(\pi_{f, W})$ for $j\neq r$ is algebraic, generated on a single geometrically connected component by the $j$-fold products of the $r$ Chern classes of the $r$ line bundles defined by the factors of automorphy attached to the non-compact archimedian places. ) Thus, to prove Theorem 3, we need only consider $H^{j}(Sh_{B, W}, \overline{\mathbf{Q}_l})(\pi_{f, W})$ where $\pi_f$ is infinite dimensional. In this case, $\pi_{\infty}$ is discrete series and $H^{j}(Sh_{B, W}, \overline{\mathbf{Q}_l})(\pi_{f, W}) \neq 0$ iff $j=r$. Let $C$ be an irreducible subquotient of $H^{j}(Sh_{B, W}, \overline{\mathbf{Q}_l})(\pi_{f, W})$. 
Then, by Proposition 6,  $C$ is a direct summand of $R_l(JL(\pi^{\prime}))$ for some 
$\pi^{\prime}$.  Note that each $ ^{[\overline{\tau_j}]} \rho^T_{l, L}$ is semisimple and 
satisfies the WMC at each finite place, since $\rho^{T}_l$ has these properties. Since, for 
$L$ as before, $R_l(JL(\pi^{\prime})|_{L} $  is a tensor product of such representations, it 
also satisfies WMC at each finite place. Hence the summand $C$ satisfies WMC and 
consequently the semisimplification of  $H^{r}(Sh_{B, W}, \overline{\mathbf{Q}_l}) $ 
as $\Gamma_L$-module also satisfies WMC at each finite place. By Proposition 3, this 
means $H^{r}(Sh_{B, W}, \overline{\mathbf{Q}_l}) $  itself satisfies WMC. \\

\subsection{Proof of Theorem 3: Langlands' Conjecture}

\noindent We recall the statement. Let $\pi^{\prime}$ be a cuspidal holomorphic 
automorphic representation of $G(\mathbf{A}_{F})$ having weight 2 at each unramified 
infinite place and with central character $|\cdot|^{-1}\Psi$ where $\Psi$ is a character of 
finite order. Let $F^{\prime}$ be the canonical field of definition of $Sh_{B,W}$. Then, 
for each finite place $v$ of $F^{\prime}$, the element

$$(\rho^{*}_{W,v}, N_{W,v})$$

\noindent of $Rep^{s}(WD_{F_v})$ defined by the restriction of the 
$\Gamma_{F^{\prime}}$ action on $H^{r}(Sh_{B, W}, 
\overline{\mathbf{Q}_l})(\pi^{\prime}_f)$ to a decomposition group for $v$, coincides 
with the class

$$m(\pi^{\prime}_f, 
W)r_B(\sigma(JL(\pi^{\prime})_{p})|_{WD_{F^{\prime}_{v}}}),$$

\noindent where 
\begin{enumerate}
\item
the non-negative integer $m(\pi^{\prime}_f, W)$ is defined by

$$ dim(H^{r}(Sh_{B, W}, \overline{\mathbf{Q}_l})(\pi^{\prime}_{f, W}))= 2^{r} 
m(\pi^{\prime}_f, W),$$

\item $p$ is the rational prime under $v$ and $JL(\pi^{\prime})$ is regarded as an 
automorphic representation of the $\mathbf{Q}$-group $Res_{F/\mathbf{Q}}(GL(2))$ 
whose Langlands parameter at $p$ is

$$\sigma(JL(\pi^{\prime}))_{p}:WD_{p}\rightarrow 
^{L}Res_{F/\mathbf{Q}}(GL(2))$$

\noindent(\cite{zeta}, 3.5), and 
\item

$$r:^{L}Res_{F/\mathbf{Q}}(GL(2))|_{F^{\prime}}\rightarrow GL(2^r, 
\mathbf{C})$$

\noindent is the representation defined by Langlands (c.f. \cite{zeta}, 5.1, 7.2).
\end{enumerate}

\noindent {\bf{Proof.}} As before, although the statement is local, for each $v$, the proof 
proceeds via the global Galois representations. In order to see clearly what is being claimed, 
we review the key definitions. \\

\noindent Let

$$R^T_l= 
r_B(Ind^{F}_{\mathbf{Q}}(\rho_{l}^{T}(JL(\pi^{\prime})))|_{\Gamma_{F^{\prime}}
}).$$

\noindent Here, for any two dimensional $l$-adic representation $\rho$ of $\Gamma_F$,  
$$Ind^{F}_{\mathbf{Q}}(\rho):\Gamma_{\mathbf{Q}}\rightarrow 
^{L}Res_{F/\mathbf{Q}}(GL(2))_l$$

\noindent is a representation of $\Gamma_{\mathbf{Q}}$ into the $l$-adic L-group 
$$^{L}Res_{F/\mathbf{Q}}(GL(2))_l.$$ This latter group is defined in general  as in 
\cite{zeta}, 3.5 using groups $\hat{G}_{l}=GL(2, \overline {\mathbf{Q}_l})$ in lieu of 
the complex groups $\hat{G}=GL(2,\mathbf{C})$. Thus, in this case,

$$^{L}Res_{F/\mathbf{Q}}(GL(2))_l= GL(2, 
\overline{\mathbf{Q}_{l}})^{Hom(F,\mathbf{R})}\times \Gamma_{\mathbf{Q}}$$

\noindent is the  semidirect product defined via the action: if $g=(g_{\tau})_{\tau \in 
Hom(F, \mathbf{R})}$ and $\eta \in \Gamma_{\mathbf{Q}}$, then $ \eta (g)_{\tau}= 
g_{\eta^{-1}\tau}$ . The homomorphism

$$I=Ind^{F}_{\mathbf{Q}}(\rho_{l}^{T}(R))$$

\noindent is defined as $$I(\eta)= ((\rho_{l}^{T}(\eta_{\overline{\tau}})_{\tau \in Hom(F, 
\mathbf{R})}), \eta)$$

\noindent where the $\overline{\tau}$ are a set of representatives in 
$\Gamma_{\mathbf{Q}}$for the $\tau$, and $\eta_{\overline{\tau}}$ is defined by the 
identity

$$ \eta\overline{\eta^{-1}\tau}=\overline{\tau}\eta_{\overline{\tau}},$$

\noindent for all $\eta$ and all $\tau$. \\

\noindent Of course, if $\eta$ fixes the Galois closure of $F$, then

$$\eta_{\overline{\tau}}=\overline{\tau}^{-1}\eta\overline{\tau},$$

\noindent so $I$ is expressed in terms of the  conjugates $^{\overline{\tau}}\rho^T_l$ of 
$\rho^T_l$.\\

\noindent We denote the  inverse image  of $\Gamma_{F^{\prime}} \subseteq 
\Gamma_{\mathbf{Q}}$ in $^{L}Res_{F/\mathbf{Q}}(GL(2))_l$ by  
$$^{L}Res_{F/\mathbf{Q}}(GL(2))_l|_{\Gamma_{F^{\prime}}}.$$

\noindent On this latter group is defined the irreducible representation $r_B$ on 
$\overline{\mathbf{Q}_{l}}^{2^r}$ . We review its construction.  Recall that  $J_{F, 
nc}= \{\tau_{1} , ..., \tau_{r}\}$ is an ordering of the set of real embeddings $J_{F, 
\mathbf{R}}\subseteq Hom(F,\mathbf{R})$ of $F$ where $B$ is split. Then on the 
connected component $GL(2, \overline{\mathbf{Q}_{l}})^{Hom(F,\mathbf{R})}$, and  
for $g=(g_{\tau})_{\tau \in Hom(F, \mathbf{R})}$,

$$r_B(g)=\otimes_{i=1}^{i=r}g_{\tau_i}.$$

\noindent  By definition of $F^{\prime}$, an $\eta \in \Gamma_{F^{\prime}}$ on $ 
Hom(F,\mathbf{R})$ defines a permutation of $J_{F, nc}$. If we define $r_B^{\eta}$ by 
the rule, $$r_B^{\eta}(g)=r_B(\eta(g)),$$
\noindent then $r_B^{\eta}$ is isomorphic to $r_B$. Let

$$P \subset  GL(2, \overline{\mathbf{Q}_{l}})^{Hom(F,\mathbf{R})}$$

\noindent be product of the groups of upper triangular matrices in each factor. Then (i) 
$\eta(P)=P$ and  (ii) $ r_B(P)$ fixes a unique line $\Lambda$ in 
$\overline{\mathbf{Q}_{l}}^{2^r}$. If $i(\eta)$ is an isomorphism satisfying, for all 
$g$,

$$i(\eta)r_B(g)= r_B^{\eta}(g)i(\eta),$$

\noindent then $i(\eta)(\Lambda)=\Lambda$. The choice of $i(\eta)$ is, by Schur's Lemma, unique up to a 
scalar, and we define  $r_B(\eta)$ to be the unique choice which leave $L$ pointwise fixed. 
The rule $r_B((g,\eta))=r_B(g)r_B(\eta)$ gives the sought representation. \\

\noindent Note that the restriction of  $R_l^T$ to the Galois closure $L$ of $F$ is just the 
representation $R_l(\pi)$ defined in Section 3. Hence $R_l^T$ is semisimple and satisfies 
WMC at each finite place $v$ of $F^{\prime}$ where $(v, l)=1$. Furthermore,  for such 
$v$, the representations $(\rho_v, N_v)$ of $WD_v$ defined by the $\Phi$-
semisimplification of the restriction of $R_l^T$ to a decomposition group 
$D_{\overline{v}}$ at $v$ coincide, since $\rho^T_l$ satisfies the Langlands 
correspondence, with

$$r_B(\sigma(JL(\pi^{\prime})_{p})|_{WD_{F^{\prime}_{v}}}).$$

\noindent Now, at this point we know that, for all $v$ prime to $l$, the representations

$$(\rho^{*}_{W,_v}, N_{W,_v})$$

\noindent of $WD_v$, defined  by the restriction to $D_{\overline{v}}$ of the  
$\Gamma_{F^{\prime}}$ module

$$H^{r}(Sh_{B, W}, \overline{\mathbf{Q}_l})(\pi^{\prime}_{f,W})$$

\noindent satisfy-whatever they may be- WMC. Thus (see 1.12), for each $v$, the nilpotent 
data $N_v$ and $N_{W,v}$ are uniquely determined by the semisimple representations 
$\rho_v$ and $\rho^{*}_{W,v}$. Hence it will suffice  to show that

$$m(\pi^{\prime}_f, W)\rho_v= \rho^{*}_{W,_v}.$$

\noindent Now, for almost all $v$, (i) $N_v=0$ and $N_{W,v}=0$, (ii) $\rho_v$ and 
$\rho^{*}_{W,v}$ are unramified, and (iii) the computation 
(\cite{ssz},\cite{zssv},\cite{brylab}) of the unramified zeta function shows exactly that this 
formula holds. Using the $l$-adic Cebotarev theorem again, we see that the semisimplified 
$\Gamma_{F^{\prime}}$-module

$$H^{r}(Sh_{B, W}, \overline{\mathbf{Q}_l})(\pi^{\prime}_{f,W})^{ss}$$

\noindent is isomorphic to   $$m(\pi^{\prime}_f, W)R_l^T.$$

\noindent Now let $v$ be any place of $F^{\prime}$ which is prime to l. Then evidently,

$$(H^{r}(Sh_{B, W}, 
\overline{\mathbf{Q}_l})(\pi^{\prime}_{f,W})|D_{\overline{v}})^{ss}$$

\noindent is isomorphic to

$$m(\pi^{\prime}_f, W)((R_l^T)|D_{\overline{v}})^{ss}.$$

\noindent Since the former gives rise to the parameter $\rho^{*}_{W, v}$ and the latter 
gives rise to $m(\pi^{\prime}_f, W)\rho_v$, we are done.


\begin{thebibliography}{999} 



\bibitem[AC]{ac}

Arthur, J.; Clozel, L., {\em {Simple algebras, base change, and the advanced theory of 
the trace formula.}} Annals of Mathematics Studies, 120. Princeton University Press, 
Princeton, NJ, 1989.



\bibitem[BR1]{brhmf}

 Blasius, D.; Rogawski, J., {\em {Motives for Hilbert modular forms.}} Invent. Math. 
114 (1993), no. 1, 55-87.



\bibitem[BR2]{zeta}

Blasius, Don; Rogawski, Jonathan D., {\em{ Zeta functions of Shimura varieties.}}  Motives (Seattle, WA, 1991),  525-571, Proc. Sympos. Pure Math., 55, Part 2, Amer. Math. Soc., Providence, RI, 1994. 



\bibitem[BrLa]{brylab} 

Brylinski, J.-L.; Labesse, J.-P., {\em{Cohomologie d'intersection et fonctions $L$ de certaines vari\'{e}t\'{e}s de Shimura.}} Ann. Sci. \'{E}cole Norm. Sup. (4)  17  (1984),  no. 3, 361--412.



\bibitem[Ca]{hmfzeta}

Carayol, Henri, {\em{Sur les repr\'{e}sentations $l$-adiques associ\'{e}es aux formes modulaires de Hilbert.}}  Ann. Sci. \'{E}cole Norm. Sup. (4)  19  (1986),  no. 3, 409--
468. 



\bibitem[C1]{annarbor}

Clozel, Laurent, {\em{ Motifs et formes automorphes: applications du principe de 
fonctorialit\'{e}.}} Automorphic forms, Shimura varieties, and $L$-functions, Vol. I (Ann Arbor, MI, 1988),  77--159, Perspect. Math., 10, Academic Press, Boston, MA, 
1990.

\bibitem[DJ]{alter}  

De Jong, A. J., {\em{Smoothness, semi-stability and alterations.}}  Inst. Hautes 
\'{E}tudes Sci. Publ. Math.  No. 83 (1996), 51--93.

\bibitem[DS]{deshalit} 

De Shalit, E.,{\em{ The p-adic monodromy-weight conjecture for p-adically uniformized varieties.}} Comp. Math. 141 (2005), 101-120.



\bibitem[D1]{dnice}

Deligne, P., {\em {Th\'{e}orie de Hodge I}}. Actes ICM Nice, t.I, pp. 425-430, Gauthier-Villars, 1970. 



\bibitem[D2]{dsv1}

Deligne, P., {\em{Travaux de Shimura}}. S\'{e}minaire Bourbaki 1970-71, Expos\'{e} 389, Lecture Notes in Mathematics, Vol. 244. Springer-Verlag, Berlin-New York, 1971.



 \bibitem[D3]{dcef}

Deligne, P., {\em{Les constantes des \'{e}quations fonctionnelles des fonctions $L$.}} 
Modular functions of one variable, II (Proc. Internat. Summer School, Univ. Antwerp, 
Antwerp, 1972),  pp. 501--597. Lecture Notes in Mathematics, Vol. 349, Springer-Verlag, Berlin, 1973.



\bibitem[D4]{weil2} 

Deligne, Pierre, {\em{La conjecture de Weil II.}} Inst. Hautes \'{E}tudes Sci. Publ. Math. No. 52 (1980), 137--252.  



\bibitem[DM]{tancats}

 Deligne, P., Milne, J., {\em {Tannakian Categories.}} Hodge cycles, motives, and Shimura varieties. Lecture Notes in Mathematics, Vol. 900. Springer-Verlag, Berlin-New York, 1982.



\bibitem[DMOS]{dmos}

 Deligne, P.; Milne, J.; Ogus, A.;  Shih, K.,  {\em {Hodge cycles, motives, and Shimura varieties.}} Lecture Notes in Mathematics, Vol. 900. Springer-Verlag, Berlin-New York, 
1982.



 \bibitem[HT]{ht}

Harris, M.; Taylor, R., {\em{The Geometry and Cohomology of Some Simple Shimura Varieties}}, Annals of Mathematics Studies, 151. Princeton University Press, Princeton, 
NJ, 2001.



\bibitem[SGA7-I]{sga7-1}

Grothendieck, A., {\em{Groupes de monodromie en g\'{e}om\'{e}trie alg\'{e}brique, I.}} Lecture Notes in Mathematics, Vol. 288, Springer Verlag, Berlin-New York, 1972.



\bibitem[Il]{illusie}

 Illusie, Luc,  {\em{Autour du th\'{e}or\'{e}me de monodromie locale.}}  P\'{e}riodes $p$-adiques (Bures-sur-Yvette, 1988).  Ast\'{e}risque  No. 223 (1994), 9-57. 



\bibitem[It1]{ito3folds} Ito, T., Weight-monodromy conjecture for certain 3-folds in 
mixed characteristic, math.NT/0212109, 2002.



\bibitem [It2]{itopadic} Ito, T., Weight-monodromy conjecture for p-adically uniformized varieties, MPI 2003-6, math.NT/ 0301202, 2003.



\bibitem[JL]{jl}

 Jacquet, H.; Langlands, R. P., {\em {Automorphic forms on ${\rm GL}(2)$.}} Lecture Notes in Mathematics, Vol. 114. Springer-Verlag, Berlin-New York, 1970. 



\bibitem[J]{jarvis}

Jarvis, Frazer, {\em{ On Galois representations associated to Hilbert modular forms.}}  J. Reine Angew. Math.  491  (1997), 199--216. 



\bibitem[Ku]{kull}

Kudla, Stephen S., {\em{The local Langlands correspondence: the non-Archimedean case.}} Motives (Seattle, WA, 1991), 365--391, Proc. Sympos. Pure Math., 55, Part 2, Amer. Math. Soc., Providence, RI, 1994. 



\bibitem[L1]{zssv} 

Langlands, R. P., {\em {On the zeta functions of some simple Shimura varieties.}} Canad. J. Math. 31 (1979), no. 6, 1121--1216.



\bibitem[L2]{lbc}

 Langlands, R.P., {\em {Base change for ${\rm GL}(2)$.}} Annals of Mathematics Studies, 96. Princeton University Press, Princeton, N.J.,1980. 

\bibitem[Oh]{ohta} 

Ohta, Masami,  {\em{On $l$-adic representations attached to automorphic forms. }} Japan. J. Math. (N.S.)  8  (1982),  no. 1, 1--47.

 \bibitem[RZ]{rz} 

Rapoport, M.; Zink, Th., {\em{Ueber die lokale Zetafunktion von Shimuravarietaeten. 
Monodromiefiltration und verschwindende Zyklen in ungleicher Charakteristik.}} Invent. 
Math.  68  (1982), no. 1, 21--101. 



 \bibitem[Ra]{rapaa} 

Rapoport, M., {\em{On the bad reduction of Shimura varieties.}}   Automorphic forms, Shimura varieties, and $L$-functions, Vol. II (Ann Arbor, MI, 1988),  253-321, Perspect. Math., 11, Academic Press, Boston, MA, 1990.



\bibitem[Re1]{ssz} 

Reimann, Harry, {\em{The semi-simple zeta function of quaternionic Shimura varieties.}} Lecture Notes in Mathematics, Vol.1657. Springer-Verlag, Berlin, 1997. viii+143 
pp.



\bibitem[Re2]{quat}  Reimann, Harry, {\em{On the zeta function of quaternionic Shimura varieties.}} Math. Ann.  317  (2000),  no. 1, 41--55. 

\bibitem[RZ]{rz}

Reimann, Harry; Zink, Thomas,  {\em {Der DieudonnŽmodul einer polarisierten abelschen Mannigfaltigkeit vom CM-Typ.}} Ann. of Math. (2)  128  (1988),  no. 3, 461--482.

\bibitem[RT]{rogtun}

 Rogawski, J. D.; Tunnell, J. B., {\em{ On Artin $L$-functions associated to Hilbert 
modular forms of weight one.}}  Invent. Math.  74  (1983),  no. 1, 1--42. 



\bibitem[Roh]{ecwd}  

Rohrlich, David E., {\em{Elliptic curves and the Weil-Deligne group.}}  Elliptic curves 
and related topics,  125--157, CRM Proc. Lecture Notes, 4, Amer. Math. Soc., 
Providence, RI, 1994. 



\bibitem[Sa]{addend}

Saito, Takeshi, {\em{Weight-monodromy conjecture for $l$-adic representations associated to modular forms.}} A supplement to: "Modular forms and $p$-adic Hodge theory" [Invent. Math.  129 (1997), no. 3, 607-620] ,  in The arithmetic and geometry of algebraic  cycles (Banff, AB, 1998), 427-431,  NATO Sci. Ser. C Math. Phys. Sci., 548, Kluwer Acad. Publ., Dordrecht, 2000. 



\bibitem[Se]{abladic}

Serre, Jean-Pierre, {\em {Abelian l-adic Representations and Elliptic Curves.}} Addison-Wesley, Redwood City, 1989.



\bibitem[ST]{serretate}

Serre, Jean-Pierre; Tate, John, {\em{Good reduction of abelian varieties.}} Ann. of Math. (2) 88 (1968) 492--517.





\bibitem[Sha]{sharam} 

 Shahidi, Freydoon, {\em{On the Ramanujan conjecture and finiteness of poles for certain $L$-functions.}} Ann. of Math. (2)  127  (1988),  no. 3, 547--584.



\bibitem[Shi]{shcurves} Shimura, Goro, {\em{ Construction of class fields and zeta functions of algebraic curves.}}  Ann. of Math. (2)  85  (1967), 58--159.



\bibitem[Tad]{tadic}

Tadi\'{c}, Marko, {\em{Classification of unitary representations in irreducible representations of general linear group (non-Archimedean case).}}  Ann. Sci. \'{E}cole 
Norm. Sup. (4)  19  (1986),  no. 3, 335--382. 



\bibitem[Ta]{ntb}

 Tate, J., {\em{Number theoretic background.}} Automorphic forms, representations and $L$-functions (Corvallis, Oregon, 1977), 3-26, Proc. Sympos. Pure Math., XXXIII, Part 2, Amer. Math. Soc., Providence, R.I., 1979. 



\bibitem[T1]{rthmf}

Taylor, R., {\em{On Galois representations associated to Hilbert modular forms.}} Invent. Math. 98 (1989), no. 2, 265-280.



\bibitem[T2]{rthmf2}

Taylor, R., {\em{On Galois representations associated to Hilbert modular forms. II.}} Elliptic curves, modular forms, \& Fermat's last theorem (Hong Kong, 1993), 185--191, Ser. Number Theory, I, Internat. Press, Cambridge, MA, 1995. 

\bibitem[TY]{tywmc}

Taylor, R.; Yoshida, Teruyoshi, {\em{ Compatibility of local and global Langlands correspondences}} preprint, math.NT/0412357, December 2004. 

\bibitem[Y]{hyhmfpds}

Yoshida, Hiroyuki, {\em{On a conjecture of Shimura concerning periods of Hilbert modular forms.}} Am. J. Math. 117 (1995), no. 4,  1019-1038. 

\bibitem[W]{ordinary}

 Wiles, A., {\em{On ordinary $\lambda$-adic representations associated to modular forms.}}  Invent. Math.  94  (1988),  no. 3, 529--573. 



\end{thebibliography}
\end{document}